\documentstyle[a4,12pt]{article}
\begin{document}
\begin{center}
{\bf Geometry and integrability of the M-LXIX  equation}
\end{center}
\begin{center}
Kuralay MYRZAKUL and R. MYRZAKULOV
\end{center}
\begin{center}
{\it Institute of Physics and Technology,
480082, Alma-Ata-82, Kazakhstan. \\
E-mail: cnlpmyra@satsun.sci.kz}
\end{center}

\section{Introduction}
The deep interrelation between many nonlinear differential equations of
the classical differential geometry of surfaces
and  modern soliton equations is well established now [1-10]. 

In this paper, we will show that the (1+1)-dimensional
spin system called the Myrzakulov LXIX (M-LXIX) equation is integrable. We do it 
by the establishing the equivalence between this and the 
Gauss-Codazzi (GC) equation. Integrability of the GC equation was established by 
V.E.Zakharov in [11].

\section{The M-LXIX equation}
Consider  the M-LXIX equation [2]
$$
{\bf S}_{t}=\frac{1}{\sqrt{S_{x}^{2}}}(-\sqrt{S_{x}^{2}-u^{2}}
{\bf S}_{x}+u{\bf S}\wedge {\bf S}_{x}) \eqno(1a)
$$
$$
u_{x}= v\sqrt{{\bf S}_{t}^{2}-u^{2}} \eqno(1b)
$$
$$
v_{t}=-{\bf S}\cdot ({\bf S}_{t}\wedge {\bf S}_{x})  \eqno(1c)
$$
where ${\bf S}=(S_{1}, S_{2}, S_{3})$ is the spin vector, 
${\bf S}^{2}=\beta=\pm 1$, $u, v$ are scalar functions. 
In this paper we study integrability of  the M-LXIX equation and its
geometry.

\section{L-equivalent of  the M-LXIX  equation}
Let us introduce the orthogonal basis defined by the vectors
$$
{\bf e}_{1}={\bf S},  \quad {\bf e}_{2}=\frac{1}{k}{\bf S}_{x},
\quad {\bf e}_{3}=\frac{1}{k}{\bf S}\wedge {\bf S}_{x}  \eqno(2)
$$
where
$$
k^{2}={\bf S}_{x}^{2}. \eqno(3)
$$
Hence and from (1) we get
$$
\left (\begin{array}{c}
e_1 \\ e_2 \\ e_3
\end{array} \right)_x= A
\left (\begin{array}{c}
e_1 \\ e_2 \\ e_3
\end{array} \right),   \quad
\left (\begin{array}{c}
e_1 \\ e_2 \\ e_3
\end{array} \right)_t= B
\left(\begin{array}{c}
e_1 \\ e_2 \\ e_3
\end{array} \right) \eqno(4)
$$
where
$$ A=
\left ( \begin{array}{ccc}
0 & k & 0 \\
-\beta k & 0 & \tau \\
0 & -\tau & 0
\end{array} \right ), \quad
B=
\left ( \begin{array}{ccc}
0 & \omega_3 & -\omega_2 \\
-\beta \omega_3 & 0 & \omega_1 \\
\beta \omega_2 & -\omega_1 & 0
\end{array} \right ). \eqno (5)
$$
In our case we put $\omega_{1}=0, \quad \beta=\pm 1$. 
The compatibility condition of these equations gives
$$
k_{t}=\omega_{3x}+\tau\omega_{2} \eqno(6a)
$$
$$
\tau_{t}=-k\omega_{2} \eqno(6b)
$$
$$
\omega_{2x}=\tau \omega_{3}. \eqno(6c)
$$
This system is the L-equivalent (Lakshmanan equivalent [1]) counterpart of
the M-LXIX equation (1).

\section{The M-LXIX equation as the particular case of the M-0
equation}
Consider the (1+1)-dimensional Myrzakulov 0 (M-0) equation [2]. This equation can be written
in the different forms. Here we work with the following form of the (1+1)-dimensional
M-0  equation
$$
{\bf e}_{1t} = \omega_{3}{\bf e}_{2} - \omega_{2} {\bf e}_{3}
\eqno(7a)
$$
$$
\tau_{t}-\omega_{1x}=
{\bf e}_{1} \cdot ({\bf e}_{1x} \wedge {\bf e}_{1t}). \eqno(7b)
$$

Let
$$
\tau=f_{x}+v, \quad \omega_{1}=f_{t} \eqno(8)
$$
where $f$ is  some function. Let $f=constant$. Then from (7) we obtain 
the equation (1). This means that the M-LXIX equation (1)
is the particular case of the M-0 equation (7).

\section{Geometry of the M-LXIX  equation}

We consider a surface immersed into the three-dimensional 
Euclidean space $R^{3}$ generated by a position vector 
${\bf r}={\bf r}(x,t)$ which is a function of two parameters -
local coordinates -  $x $ and $t$. The first and second
fundamental forms are given by
$$
I=dx^2+dy^2+dz^2=g_{ij}dx^i dx^j=
Edu^2+2Fdudv+Gdv^2     \eqno(9)
$$
and
$$
II=g_{ij}dx^{i}dx^{j}=Ldu^2+2Mdudv+Ndv^2. \eqno(10)
$$
The Gaussian curvature $K$ and  mean curvature $H$ are of the form
$$
K=\frac{b_{11}b_{22}-b_{12}^2}{g_{11}g_{22}-g_{12}^2}=
\frac{LN-M^2}{EG-F^2}.
\eqno(11)
$$
and
$$
H=\frac{1}{2}(k_1+k_2)=\frac{EN-2FM+GL}{2(EG-F^2)}.   \eqno(12)
$$
As well known that everywhere outside umbilic points the first and
second fundamental forms $I$ and $II$ can be
diagonalized simultaneously (see e.g. [10]). 
Thus choosing the curvature lines as the coordinate lines,
one, generally, has
$$
I=g_{11}dt^{2}+g_{22}dx^{2}, \quad
II=d_{11}dt^{2}+d_{22}dx^{2}. \eqno(13)
$$
In these notations the  GC equation takes the form [10]
$$
(\frac{d_{11}}{\sqrt{g_{11}}})_{x}-
\frac{d_{22}}{g_{22}}(\sqrt{g_{11}})_{x}=0
\eqno(14a)
$$
$$
(\frac{d_{22}}{\sqrt{g_{22}}})_{t}-
\frac{d_{11}}{g_{11}}(\sqrt{g_{22}})_{t}=0
\eqno(14b)
$$
$$
(\frac{(\sqrt{g_{22}})_{t}}{\sqrt{g_{11}}})_{t}+
(\frac{(\sqrt{g_{11}})_{x}}{\sqrt{g_{22}}})_{x}+
\frac{d_{11}d_{22}}{\sqrt{g_{11}g_{22}}}=0.
\eqno(14c)
$$
Following [10], we denote
$$
\psi_{1}=\frac{d_{11}}{\sqrt{g_{11}}},   \quad 
\psi_{2}=\frac{d_{22}}{\sqrt{g_{22}}},
\quad p=\frac{(\sqrt{g_{11}})_{x}}{\sqrt{g_{22}}},
\quad q=\frac{(\sqrt{g_{22}})_{x}}{\sqrt{g_{11}}}.
   \eqno(15)
$$ 
In such notations GC equation (14) takes the form
$$
\psi_{1x}=p\psi_{2} \eqno(16a)
$$
$$
\psi_{2t}=q\psi_{1} \eqno(16b)
$$
$$
q_{t}+p_{x}+\psi_{1}\psi_{2}=0. \eqno(16c)
$$
Also we have the following system [10]
$$
\tilde\psi_{1x}=p\tilde\psi_{2} \eqno(17a)
$$
$$
\tilde\psi_{2t}=q\tilde\psi_{1} \eqno(17b)
$$
where $ \tilde\psi_{1}=\sqrt{g_{11}}, \quad \tilde\psi_{2}=\sqrt{g_{22}}
$.
In such notations the three fundamental forms of the surface 
look like
$$
I=\tilde\psi_{1}^{2}dt^{2}+\tilde\psi^{2}_{2}dx^{2},\quad
II=\tilde\psi_{1}\psi_{1}dt^{2}+
\tilde\psi_{2}\psi_{2}dx^{2},\quad
III=\psi_{1}^{2}dt^{2}+\psi^{2}_{2}dx^{2}.\eqno(18)
$$

\section{Integrability of the M-LXIX equation}

Let us we compare the equations (6) and (16). We see that these equations
have  the same form after the following transformations
$$
k=\frac{\tilde\psi_{2t}}{\tilde\psi_{1}}, \quad \tau=\psi_{1},
\quad \omega_{2}=-\psi_{2},
\quad \omega_{3}=-\frac{\tilde\psi_{1x}}{\tilde\psi_{2}}. \eqno(19)
$$
Now let us rewrite the equation (16) in the form
$$
\psi_{1x}=
\frac{\tilde\psi_{1x}}{\tilde\psi_{2}}\psi_{2} \eqno(20a)
$$
$$
\psi_{2t}=\frac{\tilde\psi_{2t}}{\tilde\psi_{1}}\psi_{1} \eqno(20b)
$$
$$
(\frac{\tilde\psi_{2t}}{\tilde\psi_{1}})_{t}+
(\frac{\tilde\psi_{1x}}{\tilde\psi_{2}})_{x}+\psi_{1}\psi_{2}=0. 
\eqno(20c)
$$
Quite  recently it was shown that the GC equation (20) is integrable
by the dressing  method [11]. So as follows from these results the M-LXIX
equation (1) is  also integrable by the dressing method. 
This is the  main result of this paper. 
Finally we would like to note  that the equation (6) (and the equation
(20)) has 
the  following Lax representation
$$
\phi_x=U\phi,\quad \phi_t=V\phi \eqno(21)
$$
where
$$ U=\frac{1}{2i}\left ( \begin{array}{cc}
\tau & k \\
k & -\tau
\end{array} \right),  \quad
V=\frac{1}{2i}\left ( \begin{array}{cc}
\omega_1 & \omega_3+i\omega_2 \\
\omega_3-i\omega_2 & -\omega_1
\end{array} \right).  \eqno(22)
$$

\section{Conclusion}
Concluding, we have found the L-eqivalent counterpart (6) of the M-LXIX equation
(1). This
L-equivalent equation (6) has the same form with the GC equation (20) for the some
surface studied by Konopelchenko [10]. This means that we have idendified
to the M-LXIX equation (1) the some surface given by the GC equation (20). This
is the first main result of his paper. Integrability this GC equation is
proved by Zakharov in [11]. These results mean that  the M-LXIX equation (11)
is 
integrable by the dressing method that is our  second main result. 

\section{Acknowledgments}
This work was supported by INTAS, grant 99-1782.

\end{document}